\documentclass[12pt,a4paper]{amsart}
\usepackage{amssymb}
\usepackage{mathrsfs}

\headsep=1cm \numberwithin{equation}{section}
\newtheorem{thm}{Theorem}[section]
\newtheorem{cor}[thm]{Corollary}
\newtheorem{lem}[thm]{Lemma}
\newtheorem{prop}[thm]{Proposition}
\theoremstyle{definition}
\newtheorem{defn}[thm]{Definition}
\theoremstyle{remark}
\newtheorem{rem}[thm]{Remark}
\numberwithin{equation}{section}

\newcommand\Ass{\operatorname{Ass}}
\newcommand\mAss{\operatorname{mAss}}

\newcommand\Spec{\operatorname{Spec}}
\newcommand\Rad{\operatorname{Rad}}

\newcommand\grade{\operatorname{grade}}
\newcommand\height{\operatorname{height}}
\newcommand\qagd{\operatorname{qagd}}
\newcommand\agd{\operatorname{agd}}
\newcommand\qegd{\operatorname{qegd}}
\newcommand\qeg{\operatorname{qeg}}
\newcommand\egd{\operatorname{egd}}
\newcommand\qacogd{\operatorname{qacogd}}
\newcommand\qecogd{\operatorname{qecogd}}
\newcommand\ecogd{\operatorname{ecogd}}
\newcommand\acogd{\operatorname{acogd}}
\begin{document}\title[Quintasymptotic sequences over an ideal]
{Quintasymptotic sequences over an ideal and quintasymptotic cograde}
\author[Saeed Jahandoust and Reza Naghipour]{ Saeed Jahandoust and Reza Naghipour$^*$\\\\\\\,
\vspace*{0.5cm} Dedicated to Professor Rahim Zaare-Nahandi}
\address{Department of Mathematics, University of Tabriz, Tabriz, Iran.}
\email{saeed.e.jahan@gmail.com}
\address{Department of Mathematics, University of Tabriz, Tabriz, Iran;
and School of Mathematics, Institute for Research in Fundamental
Sciences (IPM), P.O. Box: 19395-5746, Tehran, Iran.}
\email{naghipour@ipm.ir} \email {naghipour@tabrizu.ac.ir}
\thanks{ 
This research was in part supported by a grant from IPM (No. 900130061).\\
$^*$Corresponding author: e-mail: {\it naghipour@ipm.ir} (Reza Naghipour)}%

\begin{abstract}
Let $I$ denote an ideal of a Noetherian ring $R$. The purpose of
this article is to introduce the concepts of quintasymptotic
sequences over $I$ and quintasymptotic cograde of $I$, and it is
shown that they play a role analogous to quintessential sequences
over $I$ and quintessential cograde of $I$, given in \cite{Ra1}.
Also, we show that, if $R$ is local, then the quintasymptotic
cograde of $I$ is unambiguously defined and behaves well when
passing to certain related local rings. Finally, we use this cograde
to characterize of two classes of local rings.
\end{abstract}
\maketitle
\section{Introduction}
Since the notion of regular sequences was first given, commutative
algebraists have been able to enrich their arsenal with a powerful
tool. In the early 90's several kinds of sequences have been considered
by various mathematicians working on different problems which they are
generalizations of regular sequences.
D. Katz and L. J. Ratliff, Jr., in \cite{KR} and \cite{Ra3} introduced the
interesting concepts of (essential prime divisors, essential
sequences, and essential grade of an ideal $I$ in a Noetherian ring
$R$) and (asymptotic prime divisors, asymptotic sequences, and
asymptotic grade of an ideal $I$ in a Noetherian ring $R$) and
therein they showed that these concepts are an excellent analogue
of, respectively, associated prime divisors, $R$-regular sequences,
and the standard grade of $I$, in the classical theory.

On the other hand, D. Rees introduced the important concept of an
asymptotic sequence over an ideal $I$ in a Noetherian local ring $R$
in \cite{Re}, and in \cite{MR1}, S. McAdam and L. J. Ratliff, Jr.,
showed that asymptotic sequences over an ideal $I$ in a Noetherian
ring $R$ and asymptotic cograde of $I$ (when $R$ is local) have same
useful properties, and several bounds on this cograde were
established in \cite{MR1}. The main purpose of the present article
is to introduce the concepts of quintasymptotic sequences over an
ideal $I$ in a Noetherian ring $R$ and quintasymptotic cograde of
$I$. We show that, quintasymptotic sequences over $I$ behave nicely
when passing to certain rings related to $R$ and that the
quintasymptotic cograde of $I$ is well defined (when $R$ is local)
and satisfied certain rather natural inequalities. Also, we show
that if $R$ is local, then any two maximal quintasymptotic sequences
over $I$ have the same length, and
\begin{center}
$\qacogd(I)=\mathrm{min}\{\dim R^\ast/{IR^\ast +z}\mid
z\,\,\text{is a minimal prime of}\,\, R^\ast\}.$
\end{center}
Finally, we show that, for every ideal $I$ in a complete local ring
$R$, $\agd(I)+\qacogd(I)=\dim R$ if and only if, for every prime
ideal $\frak p$ of $R$ with $\dim R/\frak p=1$, $\agd(\frak
p)+\qacogd(\frak p)=\dim R$, if and only if $R$ has a unique minimal
prime divisor of zero.

Throughout this paper, all rings considered will be commutative and
Noetherian and will have non-zero identity elements. Such a ring
will be denoted by $R$, and the terminology is, in general, the same
as that in \cite{BH}, \cite{Ma} and \cite{Nag1}.

Let $I$ be an ideal of $R$. We denote by $\mathscr{R}$ the {\sl Rees
ring} $R[u, It] :=\oplus_{n\in \mathbb{Z}}I^nt^n$ of $R$ w.r.t. $I$,
where $t$ is an indeterminate and $u =t^{-1}$. Also, the {\it
radical} of $I$, denoted by $\Rad(I)$, is defined to be the set
$\{x\in R: x^n \in I$ for some $n\in\mathbb{N}\}$. For each
$R$-module $L$, we denote by $\mAss_RL$ the set of minimal primes of
$\Ass_RL$. If $(R,\frak{m})$ is Noetherian local, then $R^*$ denotes
the completion of $R$ with respect to the $\frak{m}$-adic topology.
Then $R$ is said to be an {\it unmixed} (resp. {\it quasi-unmixed})
ring if for every $\frak{p}\in \Ass_{R^*}R^*$ (resp. $\frak{p}\in
\mAss_{R^*}R^*$), the condition $\dim R^*/\frak{p}=\dim R$ is
satisfied. More generally, if $R$ is not necessarily local, $R$ is a
{\it locally unmixed} (resp. {\it locally quasi-unmixed}) ring if
for any $\frak{p} \in \Spec R$, $R_\frak{p}$ is an unmixed (resp.
quasi-unmixed) ring. A prime ideal $\frak{p}$ of $R$ is called a
{\it quintessential} (resp. {\it quintasymptotic})  prime ideal of
$I$ precisely when there exists $\frak{q}\in
\Ass_{R^*_\frak{p}}R^*_\frak{p}$ (resp. $\frak{q}\in
\mAss_{R^*_\frak{p}}R^*_\frak{p}$) such that $\Rad(IR^*_\frak{p}+
\frak{q})= \frak{p}R^*_\frak{p}$. The set of {\it quintessential}
(resp. {\it quintasymptotic}) primes of $I$ is denoted by $Q(I)$
(resp. $\overline{Q^*}(I)$).
 Then the {\it essential} (resp. {\it asymptotic}) primes of $I$, denoted by $E(I)$ (resp.  $\overline{A^*}(I)$),
is defined to be the set $\{\frak{q}\cap R\mid\, \frak{q}\in Q(u\mathscr{R})\}$
(resp. $\{\frak{q}\cap R\mid\, \frak{q}\in \overline{Q^*}(u\mathscr{R})\}$).
Finally, we shall use $A^*(I)$ to denote the ultimately constant value of the sequence
$\Ass_RR/I^n$, which is well defined finite set of prime ideals, (cf. \cite{Br1}).

 A brief summary of the contents of this paper will now be given. Let $R$ be a commutative ring
and $I$ an arbitrary ideal of $R$.  In the second section, the
notion of the quintasymptotic sequences over $I$ is introduced, and
it is shown that most of the basic properties of the quintessential
sequences (resp.  essential sequences) over $I$, given in \cite{Ra1}
(resp. \cite{KR}), extend to the quintasymptotic sequences over $I$.
In fact, it is shown in this section that the quintasymptotic
sequence over $I$ behaves nicely with respect to passing to certain
rings related to $R$. In the third section, the concept of
quintasymptotic cograde of an ideal developed and it is shown that
most of bounds on the essential, asymptotic and quintessential
cograde of $I$ given in \cite{KR}, \cite{MR1} and \cite{Ra1} have a
valid analogous for the quintasymptotic cograde of $I$. Finally, in
section 4 we characterize of two class of local rings by using the
quintasymptotic cograde.


\section{Quintasymptotic sequences over an ideal}

In this section we introduce the notion of the quintasymptotic
sequences over an ideal $I$ of a Noetherian ring $R$ and show that
they have most of basic properties enjoyed by quintessential
sequences, essential sequences and  asymptotic sequences over $I$.
We begin with the following definitions.

\begin{defn}
Let $I$ and $\frak p$ be ideals of a Noetherian ring $R$ such that $\frak p$ is prime.
Then $\frak p$ is called a {\it quintasymptotic} (resp. {\it quintessential}) {\it prime ideal of} $I$
precisely when there exists $z\in \mAss_{R^*_\frak{p}}R^*_\frak{p}$ (resp. $z\in \Ass_{R^*_\frak{p}}R^*_\frak{p}$)
 such that
$\Rad(IR^*_\frak{p}+ z)=\frak{p}R^*_\frak{p}$.  The set of
quitasymptotic (resp. quintessential) primes of $I$ is denoted by
$\overline{Q^*}(I)$ (resp. $Q(I)$).
\end{defn}

\begin{defn}
Let $I$ denote an ideal of a Noetherian ring.  A sequence
$\mathbf{x}= x_1, \ldots, x_n$ of elements of $R$ is called a {\it
quintasymptotic} (resp. {\it quintessential}) {\it sequence over}
$I$ if, $(I,(\mathbf{x})) \neq R$ and for all $1\leq i\leq n$, we
have $x_i\not\in \bigcup\{\frak{p}\in \overline{Q^*}((I,(x_1,
\ldots, x_{i-1})))$ (resp. $x_i\not\in \bigcup\{\frak{p}\in
Q((I,(x_1, \ldots, x_{i-1})))$. A quintasymptotic (resp.
quintessential) sequence over $(0)$ is simply called a {\it
quintasymptotic} (resp. {\it quintessential}) {\it sequence in} $R$
\end{defn}
A quintasymptotic (resp. quintessential) sequence $\mathbf{x}= x_1,
\ldots, x_n$ of elements of $R$ over $I$ is {\it maximal} if $x_1,
\ldots, x_n, x_{n+1}$ is not a quintasymptotic (resp.
quintessential) sequence over $I$ for any $x_{n+1}\in R$. If $R$ is
local, then it is shown in Theorem 3.2 (resp. {\cite[Theorem
3.2]{Ra1}}) that all maximal quintasymptotic (resp. quintessential)
sequences over $I$ have the same length. This allows us to introduce
the fundamental notion of {\it quintasymptotic cograde} (resp. {\it
quintessential cograde}), $\qacogd(I)$ (resp. $\qecogd(I)$), of $I$.
Also, it is shown in Corollary 2.15 (resp. {\cite[Proposition
4.3]{MR2}}) that all maximal quintasymptotic (resp. quintessential)
sequences with respect to coming from $I$ have the same length.
Therefore, we define the fundamental notion of {\it quintasymptotic
grade} (resp. {\it quintessential grade}),
$\qagd(I)$ (resp. $\qegd(I)$) of $I$.\\

The following lemma is needed in the proof of the main results of this paper.
\begin{lem}
Let $I$ and $J$ be ideals in a Noetherian ring $R$. Then the
following hold:

{\rm(i)} If $\frak p$ is a minimal prime divisor of $I$, then $\frak
p\in \overline{Q^\ast}(I)$.

{\rm(ii)} If $\Rad (I)=\Rad (J)$, then
$\overline{Q^\ast}(I)=\overline{Q^\ast}(J)$.

{\rm(iii)} $\overline{Q^\ast}(I)\subseteq \overline{A^\ast}(I)\cap
Q(I)$ and $\overline{A^\ast}(I)\cup Q(I)\subseteq E(I)\subseteq
A^\ast (I)$.

{\rm(iv)} If $I\subseteq \frak p \in \rm{Spec}$ $R$ and $S$ is a
multiplicatively closed subset of $R$ which is disjoint from $\frak
p$, then $\frak p\in \overline{Q^\ast}(I)$ if and only if $\frak
pR_S\in \overline{Q^\ast}(IR_S)$.

{\rm(v)} $\frak p\in \overline{Q^\ast}(I)$ if and only if there is
$z\in \mAss_RR$ such that  $z\subseteq \frak p$ and $\frak p/z \in
\overline{Q^\ast}(I(R/z))$.

{\rm(vi)} If $z\in \mAss_RR$ and $\frak p$ is a minimal over $I+z$,
then $\frak p\in \overline{Q^\ast}(I)$.

{\rm(vii)} Let the ring $T$ be a faithfully flat Noetherian
extension of $R$. Let $\frak q$ be a prime ideal in $T$ such that
$\frak p=\frak q \cap R$. If $\frak q\in \overline{Q^\ast}(IT)$,
then $\frak p\in \overline{Q^\ast}(I)$ and $\frak q\in
\overline{Q^\ast}(\frak pT)$. Moreover, if $\frak p\in
\overline{Q^\ast}(I)$ and $\frak q$ is minimal over $\frak pT$, then
$\frak q\in \overline{Q^\ast}(IT)$.

{\rm(viii)} Let the ring $T$ be a finite module extension of $R$. If
$\frak p\in \overline{Q^\ast}(I)$, then there is a $\frak q\in
\overline{Q^\ast}(IT)$ such that $\frak q\cap R=\frak p$. Moreover, if  all
minimal primes in $T$ lies over a minimal prime in $R$, then the
converse holds.

{\rm(ix)} $\overline{Q^\ast}((I,X)R[X])=\{(\frak p,X)R[X]\mid \frak
p\in \overline{Q^\ast}(I)\}$.
\end{lem}
 \proof (i) and (ii) follow readily from definition. (iii)-(viii)
are proved in {\cite[Lemmas 2.1, 3.4 and Propositions 3.6,
3.8]{McA3}}. To prove (ix), let $\frak p \in\overline{Q^\ast}(I)$.
Then by (iv),(v) and (vii), we may assume that $R$ is a complete
local domain with maximal ideal $\frak p$. Then, in view of
{\cite[Propositions 6 and 7]{Nag2}}, $R$ and $R[X]$ are locally
unmixed. Thus $\overline{Q^\ast}(J)=Q(J)$ for all ideals $J$ in $R$
and in $R[X]$. Therefore $(\frak p,X)R[X]\in
\overline{Q^\ast}((I,X)R[X])$, by {\cite[Lemma
2.7]{Ra1}}. The other inclusion is similar. \qed\\

The next result is a consequence of  Lemma 2.3(i) and Definition
2.2.
\begin{cor}
Let $I$ be an ideal in a Noetherian ring $R$ and let ${\bf x}=x_1, \dots, x_n$ be a
sequence of elements of $R$.

{\rm(i)} If ${\bf x}$ is a quintasymptotic sequence over $I$, then
$\height(I,({\bf x}))\geq\height I +n$. Therefore by the
\emph{Generalized Principal Ideal Theorem} if  ${\bf x}$ is a
quintasymptotic sequence in $R$, then $\height(({\bf x}))=n$.

{\rm(ii)} The  sequence  ${\bf x}$ is  a  maximal quintasymptotic
sequence over $I$ if and only if ${\bf x}$ is a quintasymptotic
sequence over $I$ and for each maximal ideal $\frak m$ in $R$
containing $(I,({\bf x}))$ it holds that $\frak m\in
\overline{Q^\ast}((I, ({\bf x})))$.
\end{cor}

The following proposition shows that the quintasymptotic sequences over
an ideal are well behavior when passing to localization.

\begin{prop}
Let $I$ be an ideal in a Noetherian ring $R$ and let ${\bf x}=x_1, \dots, x_n$  be
a sequence of elements of $R$. Then the following statements hold:

{\rm(i)} If $\mathbf{x}$ is a quintasymptotic sequence over $I$ and
$S$  a multiplicatively closed subset of $R$ such that
$(I,(\mathbf{x}))R_S\neq R_S$, then the image of $\mathbf{x}$ in
$R_S$ is a quintasymptotic sequence over $IR_S$. The converse holds
if for all  $\frak p\in \bigcup \{\frak q\in
\overline{Q^\ast}((I,(x_1,\dots, x_{i})));i=0,\dots, n-1\}$, we have
$\frak pR_S\neq R_S$.

{\rm(ii)} If $\mathbf{x}$ is a maximal quintasymptotic sequence over
$I$, then for each maximal ideal $\frak m$ in $R$ containing
$(I,(\mathbf{x}))$ it holds that the image of $\mathbf{x}$ in
$R_\frak m$ is a maximal quintasymptotic sequence over $IR_\frak m$.
The converse holds if the $x_i$ are contained in the Jacobson
radical of $R$.
\end{prop}

\proof  (i) follows from Lemma 2.3(iv). The first statement in (ii)
follows from part (i) and Corollary 2.4(ii). For the last statement
in (ii) it will first be shown that \textbf{x} is a quintasymptotic
sequence over $I$. To this end, suppose the contrary is true. Then
there exists $i$ such that $x_i\in \frak p$ for some $\frak p \in
\overline{Q^\ast}((I,(x_1,\dots,x_{i-1})))$. Let $\frak m$ be a
maximal ideal in $R$ containing $\frak p$. Then the hypotheses
implies that $(I,(\textbf{x}))R\subseteq \frak m$, and so the
hypotheses  and Lemma 2.3(iv) imply that the image of $x_i$ is in
$\frak pR_\frak m\in \overline{Q^\ast}((I,(x_1,...,x_{i-1}))R_\frak
m)$. But this implies that the image of \textbf{x} in $R_\frak m$ is
not quintasymptotic sequence over $IR_\frak m$, in contradiction to
the hypotheses. Therefore $\textbf{x}$ is a quintasymptotic sequence
over $I$. If $\frak m$ is a maximal ideal in $R$ containing
$(I,(\textbf{x}))$, then by Corollary 2.4(ii) and Lemma 2.3(iv) we
have $\frak m\in \overline{Q^\ast}((I,(\textbf{x})))$, and so by
Corollary 2.4(ii),
\textbf{x} is a maximal quintasymptotic sequence over $I$.\qed \\

The following result shows that the quintasymptotic sequences over
an ideal $I$ of a Noetherian ring $R$ are well behavior when passing to the factor rings modulo
minimal primes of $R$.

\begin{prop}
Let $I$ be an ideal in a Noetherian ring $R$ and let ${\bf x}=x_1, \dots, x_n$ be
a sequence of elements of  $R$. Then the following statements hold:

{\rm(i)} ${\bf x}$ is a quintasymptotic sequence over $I$ if and
only if the image of ${\bf x}$ in $R/z$ is a quintasymptotic
sequence over $I(R/z)$ for all $z\in \mAss_RR$.

{\rm(ii)} ${\bf x}$ is a maximal quintasymptotic sequence over $I$
if and only if the image of ${\bf x}$ in $R/z$ is a quintasymptotic
sequence over $I(R/z)$ for all $z\in \mAss_RR$ and for all maximal
ideals $\frak m$ in $R$ containing $(I,({\bf x}))$, there exists
$z\in \mAss_RR$ such that $z\subseteq \frak m$ and
 $\frak m/z\in \overline{Q^\ast}((I,({\bf x}))(R/z))$.

\proof It follows readily from Lemma 2.3(v) and Corollary 2.4(ii). \qed \\
\end{prop}

The next result shows that the quintasymptotic sequences over an
ideal are well behavior when passing to faithfully flat Noetherian
extension rings of $R$. Before bringing of it, let us recall the following
definition.

\begin{defn}
Let $R\subseteq T$ be Noetherian rings.

{\rm (i)} We say that $R$ is {\it dominated by} $T$ if, for every proper ideal $I$ of $R$, we have $IT\neq T$ and
every maximal ideal of $T$ lies over a maximal ideal of $R$.

{\rm (ii)} We say that the {\it Theorem of Transition} holds for rings $R$ and $T$ if,  $R$ is dominated by $T$ and
if $\frak q$ is a primary ideal of $R$ such that $\Rad(\frak q)$ is a maximal ideal,
say $\frak m$, then ${\rm length}_T T/{\frak q T}$ is finite and that
$${\rm length}_T T/{\frak qT}=({\rm length}_T T/{\frak mT})({\rm length}_R R/{\frak q}).$$
\end{defn}

\begin{prop}
 Let $R\subseteq T$ be a faithfully flat extension
of Noetherian rings. Let $I$ be an ideal of $R$ and let
${\bf x}=x_1, \dots, x_n$ be a sequence of elements of $R$. Then the following hold:

{\rm(i)} ${\bf x}$ is a quintasymptotic sequence over $I$ if and
only if ${\bf x}$ is a quintasymptotic sequence over $IT$.

{\rm(ii)} If $R\subseteq T$ satisfy the Theorem of Transition, then
${\bf x}$ is a maximal quintasymptotic sequence over $I$ if and
only if ${\bf x}$ is a maximal quintasymptotic sequence over
$IT$.
\end{prop}
\proof (i) follows immediately from Lemma 2.3(vii). In order to
prove (ii), let $\mathbf{x}$ be  maximal quintasymptotic sequence
over $I$ and $\frak n$  a maximal ideal of $T$ containing
$(I,(\mathbf{x}))T$. Then by (i), $\mathbf{x}$ is a quintasymptotic
sequence over $IT$. Let $\frak m:=\frak n \cap R$.  Since $R$ is
dominated by $T$, it follows that $\frak m$ is a maximal ideal
containing $(I,(\mathbf{x}))$, and so $\frak m\in
\overline{Q^\ast}((I,(\mathbf{x})))$, by  Corollary 2.4(ii).
Therefore $\frak n\in \overline{Q^\ast}((I,(\mathbf{x}))T)$ by Lemma
2.3(vii) and hence, $\mathbf{x}$ is a maximal quintasymptotic
sequence over $IT$, by Corollary 2.4(ii).

Now, let $\mathbf{x}$ be  maximal quintasymptotic sequence over
$IT$. Then by (i), $\mathbf{x}$ is a quintasymptotic sequence over
$I$. Let $\frak m$ be a maximal ideal of $R$ containing
$(I,(\mathbf{x}))$ and let $\frak n$ be a maximal ideal in $T$
containing $\frak mT$. Then $\frak n\in
\overline{Q^\ast}((I,(\mathbf{x}))T)$ by  Corollary 2.4(ii). In
other hand $\frak m=\frak n\cap R$,  and so $\frak m\in
\overline{Q^\ast}((I,(\mathbf{x})))$ by Lemma 2.3(vii). Therefore
$\mathbf{x}$ is a maximal quintasymptotic
sequence over $I$, by Corollary 2.4(ii). \qed \\

The next result shows that the quintasymptotic sequences over an
ideal are well behavior when passing to finite extension rings of $R$.

\begin{prop}
Let $R\subseteq T$ be Noetherian rings, with $T$ a finitely generated $R$-module. Let  $I$
be an ideal of $R$ and let  ${\bf x}=x_1, \dots, x_n$ be a sequence of elements of $R$. Then
the following statements hold:

{\rm (i)} If $\mathbf{x}$ is a quintasymptotic sequence over $IT$, then
$\mathbf{x}$ is a quintasymptotic sequence over $I$.

{\rm (ii)} If every minimal prime of $T$ lies over a minimal prime
in $R$, then $\mathbf{x}$ is a quintasymptotic sequence over $I$ if
and only if $\mathbf{x}$ is a quintasymptotic sequence over $IT$.

{\rm (iii)} If every minimal prime of $T$ lies over a minimal prime
in $R$, then $\mathbf{x}$ is a maximal quintasymptotic sequence over
$I$ if and only if $\mathbf{x}$ is a quintasymptotic sequence over
$IT$ and for each maximal ideal $\frak m$ in $R$ that contains
$(I,(\mathbf{x}))$, there exists a prime ideal $\frak n$ in $T$ such
that $\frak m=\frak n\cap R$ and $\frak n\in
\overline{Q^\ast}((I,(\mathbf{x}))T)$.
\end{prop}
\proof It follows readily from Lemma 2.3 (viii) and Corollary 2.4(ii). \qed \\

The next proposition is concerned with the quintasymptotic sequences over
$I$ and $IR[X]$, where $X$ is an indeterminate over $R$.
\begin{prop}
Let I be an ideal in a Noetherian ring R and let ${\bf x}=x_1, \dots, x_n$ be a
sequence of elements of R. Then the following statements are equivalent:

{\rm (i)} The sequence $\mathbf{x}$ is a (resp. maximal) quintasymptotic sequence over I.

{\rm (ii)} The sequence  $x_1, \dots, x_i, X, x_{i+1}, \dots, x_n$ is a (resp. maximal) quintasymptotic
sequence over $IR[X]$ for some $i=0,1,\dots, n$.

{\rm (iii)} The sequence $x_1, \dots, x_i, X, x_{i+1}, \dots, x_n$ is a (resp. maximal) quintasymptotic
sequence over $IR[X]$ for every $i=0,1, \dots, n$.
\end{prop}

\proof  In view of Lemma 2.3(vii),  for $j=0,1,\dots, i$, we have
$$\overline{Q^\ast}((I, (x_1, \dots, x_j))R[X])=\{\frak pR[X]\mid \frak p\in
\overline{Q^\ast}((I, (x_1, \dots, x_j))\},$$ (note that, for an
ideal $J$ in $R$, the prime divisors of $JR[X]$ are  $\frak pR[X]$
such that $\frak p$ is  a prime divisor of $J$). Also, it is clear
that $X$ is not in any prime divisor of $(I, (x_1, \dots,
x_i))R[X]$. Moreover, for $k =0,1, \dots, n-i$, we have
$$\overline{Q^\ast}((I, (x_1, \dots, x_i,X,x_{i+1}, \dots, x_{i+k}))R[X])=\{(\frak
p,X)R[X]\mid \frak p\in \overline{Q^\ast}((I,(x_1, \dots,
x_{i+k})))\},$$ by Lemma 2.3(ix). Now,  the result follows from this
and Corollary 2.4(ii). Note that  the maximal ideals of $R[X]$
containing $(I,X)R[X]$ are the ideals $(\frak m,X)R[X]$ such that
$\frak m$ is a maximal ideal of $R$ containing
$I$.\qed\\

The following result is concerned with quintasymptotic sequences
over ideals with the same radical.

\begin{prop}
Let $I$ and $J$ be ideals in a Noetherian ring $R$ such that $\Rad(I)=\Rad(J)$,
and let ${\bf x}=x_1, \dots, x_n$  be a sequence of elements of $R$. Then the
following statements hold:

{\rm (i)} $\mathbf{x}$ is a quintasymptotic sequence over $I$ if and
only if $\mathbf{x}$ is a quintasymptotic sequence over $J$.

{\rm (ii)} $\mathbf{x}$ is a maximal quintasymptotic sequence over
$I$ if and only if $\mathbf{x}$ is a maximal quintasymptotic
sequence over $J$.
\end{prop}
\proof As $\Rad (I)=\Rad (J)$, it follows that
$$\Rad ((I, (x_1,\dots, x_i)))= \Rad ((J, (x_1, \dots, x_i))),$$ for all
$i=1,\dots, n$,  and in view of Lemma 2.3(ii), we have
$$\overline{Q^\ast}((I, (x_1, \dots, x_i)))=\overline{Q^\ast}((J, (x_1, \dots, x_i))),$$
for all $i=1, \dots, n$. Therefore (i) is true by Definition 2.2.

 Now, let $\mathbf{x}$ be a maximal quintasymptotic sequence over $I$.  Then by
(i), $\mathbf{x}$ is a quintasymptotic sequence over $J$.  Let
$\frak m$ be a maximal ideal of $R$ containing $(J,(\mathbf{x}))$.
Then, $(I,(\mathbf{x}))\subseteq \frak m$. Since $\mathbf{x}$ is a
maximal quintasymptotic sequence over $I$, it follows from Corollary
2.4(ii) that $\frak m\in \overline{Q^\ast}((I,(\mathbf{x})))$.
Therefore $\frak m\in \overline{Q^\ast}((J,(\mathbf{x})))$, and so
by Corollary 2.4(ii), $\mathbf{x}$ is a maximal quintasymptotic
sequence over $J$.
The converse will be proved similarly.\qed \\

The next remark, which gives us some additional basic informations
concerning quintasymptotic sequences over an ideal, follows
immediately from  the Definition 2.2 and Lemma 2.3(ii).

\begin{rem}
Let $I$ be an ideal in a Noetherian ring $R$ and let $\textbf{x}=x_1, \dots, x_n$ be
a sequence of elements of $R$. Then the following statements are equivalent:

(i) $x_1,\dots, x_n$ is a quintasymptotic sequence over $I$.

(ii) ${x_1}^{m_1},\dots, {x_n}^{m_n}$ is a quintasymptotic sequence
over $I$ for some positive integers $m_i$.

(iii) ${x_1}^{m_1},\dots,{x_n}^{m_n}$ is a quintasymptotic sequence
over $I$ for all positive integers $m_i$.

(iv) There exists an integer $i\in\{0,\dots, n-1\}$ such that $x_1,\dots, x_i$ is a
quintasymptotic sequence over $I$ and $x_{i+1},\dots, x_n$ is
quintasymptotic sequence over $(I, (x_1,\dots, x_i))$.

(v) For all $i\in \{0,\dots, n-1\}$, $x_1,\dots, x_i$ is a quintasymptotic
sequence over $I$ and $x_{i+1},\dots, x_n$ is a quintasymptotic sequence
over $(I,(x_1,\dots, x_i))$.
\end{rem}

In the remainder of this section, we examine the quintasymptotic
sequences over zero ideal. Before bringing the next results we
recall a sequence $\mathbf{x}= x_1, \ldots, x_n$ of elements of $R$
is called an {\it asymptotic} (resp. {\it essential}) {\it sequence
over} $I$ if, $(I,(\mathbf{x})) \neq R$ and for all $1\leq i\leq n$,
we have $x_i\not\in \bigcup\{\frak{p}\in \overline{A^*}((I,(x_1,
\ldots, x_{i-1})))$ (resp. $x_i\not\in \bigcup\{\frak{p}\in
E((I,(x_1, \ldots, x_{i-1})))$. An asymptotic (resp. essential)
sequence over $(0)$ is simply called an {\it asymptotic} (resp. {\it
essential}) {\it sequence in} $R$. An asymptotic (resp. essential)
sequence $\mathbf{x}= x_1, \ldots, x_n$ of elements of $R$ over $I$
is {\it maximal} if $x_1, \ldots, x_n, x_{n+1}$ is not an asymptotic
(resp. essential) sequence over $I$ for any $x_{n+1}\in R$. If $R$
is local, then it is shown in {\cite[Theorem 1.9]{K}} (resp.
{\cite[Theorem 4.1]{KR}}) that all maximal asymptotic (resp.
essential) sequences over $I$ have the same length. This allows us
to introduce the fundamental notion of {\it asymptotic cograde}
(resp. {\it essential cograde}), $\acogd(I)$ (resp. $\ecogd(I)$), of
$I$. Also, it is shown in {\cite[Theorem 3.1]{Ra3}} (resp.
{\cite[Proposition 3.10]{KR}}) that all maximal asymptotic (resp.
essential) sequences with respect to coming from $I$ have the same
length. Therefore, we define the fundamental notion of {\it
asymptotic grade} (resp. {\it essential
grade}), $\agd(I)$ (resp. $\egd(I)$) of $I$.\\

The next proposition is a consequence of {\cite[Theorem 3.1] {KR}}
and {\cite[Proposition 4.1] {Ra1}}.

\begin{prop}
Let $R$ be a Noetherian ring such that $\Ass_{R_{\frak p}^\ast}
R_{\frak p}^\ast$ has no embedded prime ideals for all $\frak p\in
\Spec R$. Let $I$ be an ideal of $R$ and let $\mathbf{x}=x_1, \dots,
x_n$ be a sequence of elements of $R$. Then
$\overline{Q^\ast}(I)=Q(I)\subseteq \overline{A^\ast}(I)=E(I)$.
Moreover, the following statements are equivalent:

$\rm (i)$ $\mathbf{x}$ is a quintasymptotic sequence in $R$;

$\rm (ii)$ $\mathbf{x}$ is an asymptotic sequence in $R$;

$\rm (iii)$ $\mathbf{x}$ is a quintessential sequence in $R$;

$\rm (iv)$ $\mathbf{x}$ is an essential sequence in $R$.
\end{prop}

\proof It is easy to see that $\overline{Q^\ast}(I)=Q(I)$. Also, by
{\cite[Theorem 3.1] {KR}} we have $Q(I)\subseteq
\overline{A^\ast}(I)=E(I)$. Now, in view of {\cite[Proposition 4.1]
{Ra1}} and {\cite[Proposition 3.10] {KR}}, we have $\rm
(iii)\Longleftrightarrow \rm (iv)$ and $\rm
(ii)\Longleftrightarrow\rm (iii)$, and the first statement shows that $\rm (i)\Longleftrightarrow\rm (iii)$. \qed \\

Now we show that the quintasymptotic grade of $I$ is well defined and
equals with asymptotic grade of $I$.

\begin{prop}
Let $I$ be an ideal of a Noetherian ring $R$. If $I$ is generated by
a quintasymptotic sequence of elements of  $R$, then
$\overline{Q^\ast}(I)=\overline{A^\ast}(I)$.
\end{prop}

\proof In view of Lemma 2.3(iii) it suffices to us show that
$\overline{A^\ast}(I)\subseteq \overline{Q^\ast}(I)$. To do this,
let $\mathbf{x}=x_1, \dots, x_n$ be a quintasymptotic sequence of
elements of  $R$ and $I=(\textbf{x})$. Let $\frak p\in
\overline{A^\ast}(I)$. Then $\frak pR_\frak p\in
\overline{A^\ast}(IR_\frak p)$. Further, if $\frak pR_\frak p\in
\overline{Q^\ast}(IR_\frak p)$, then by Lemma 2.3(iv), $\frak p\in
\overline{Q^\ast}(I)$. Also, $IR_\frak p$ is generated by a
quintasymptotic sequence in $R_\frak p$, by Proposition 2.5. Thus we
may assume that $R$ is a local ring with maximal ideal $\frak p$.
Moreover $\frak pR^\ast\in \overline{A^\ast}(IR^\ast)$ by
{\cite[Proposition 3.18] {McA1}}, and if $\frak pR^\ast\in
\overline{Q^\ast}(IR^\ast)$, then $\frak p\in \overline{Q^\ast}(I)$,
by Lemma 2.3(vii). Also, $IR^\ast$ is generated by a quintasymptotic
sequence in $R^\ast$ by Proposition 2.8. Therefore, we may assume
that $(R,\frak p)$ is a  complete local ring. Finally, by
{\cite[Proposition 3.18] {McA1}} there exists $z\in \mAss_RR$ such
that $\frak p/z \in \overline{A^\ast}(I(R/z))$,  and if $\frak p/z
\in \overline{Q^\ast}(I(R/z))$, then $\frak p\in
\overline{Q^\ast}(I)$, by Lemma 2.3(v). Also, $I(R/z)$ is generated
by a quintasymptotic sequence in $R/z$ by Proposition 2.6.
Consequently, we may assume that $(R,\frak p)$ is a complete local
domain. Now, as by {\cite[Proposition 6] {Nag2}} $R$ is a locally
unmixed, it follows from Proposition 2.13 that $I$ is generated by a
quintessential sequence in $R$ and
 $\overline{Q^\ast}(I)=Q(I)\subseteq
\overline{A^\ast}(I)=E(I)$. Whence  by {\cite[Theorem 2.5]{KR}}, we have
$\overline{Q^\ast}(I)=\overline{A^\ast}(I)$, as required.\qed \\

\begin{cor}
Let $R$ be a Noetherian ring and let $\mathbf{x}=x_1, \dots, x_n$ be a sequence of elements of $R$.
Then $\mathbf{x}$ is an asymptotic sequence in $R$ if and only if
$\mathbf{x}$ is a quintasymptotic sequence in $R$. In particular
${\rm agd}(I)={\rm qagd}(I)$ for all ideals $I$ in $R$.
\end{cor}
\proof It follows from Lemma 2.3(iii) and Proposition 2.14.\qed \\

\begin{cor}
Let $R$ be a Noetherian ring such that $\Ass_{R_{\frak m}^\ast}
R_{\frak m}^\ast$ has no embedded prime ideals for all maximal
ideals $\frak m$ in $R$. Let $I$ be an ideal of $R$ and let
$\mathbf{x}=x_1, \dots, x_n$ be a sequence of elements of $R$. Then
the following statements are equivalent:

$\rm (i)$ $\mathbf{x}$ is a quintasymptotic sequence in $R$.

$\rm (ii)$ $\mathbf{x}$ is an asymptotic sequence in $R$.

$\rm (iii)$ $\mathbf{x}$ is a quintessential sequence in $R$

$\rm (iv)$ $\mathbf{x}$ is an essential sequence in $R$.

In particular, $\qagd(I)=\agd(I)=\qegd(I)=\egd(I)$.
\end{cor}
\proof It follows from Corollary 2.15, {\cite[Proposition 4.1] {Ra1}}
and {\cite[Proposition 3.10] {KR}}.\qed

\begin{rem}
It has been shown in {\cite[Lemma 6.13] {McA1}} that, if $\frak p\in
\overline{A^\ast}(I)$, $\mathbf{x}=x_1, \dots, x_n$ is an asymptotic
sequence over $I$ and $\frak q$ a minimal prime divisor of $(\frak
p,(\mathbf{x}))$, then $\frak q\in
\overline{A^\ast}((I,(\mathbf{x})))$. Also, this statement is proved
in {\cite[Theorem 5.1] {KR}} for essential primes and essential
sequences. But it is not true for quintasymptotic primes and
quintasymptotic sequences, because in otherwise it is not hard to
show that it holds for quintessential primes and quintessential
sequences, and this contradicts {\cite[Example 7.3] {KR}}.
\end{rem}

The next corollary is weaker result than above remark for
quintasymptotic primes and quintasymptotic sequences.

\begin{cor}
Let $\mathbf{x}=x_1, \dots, x_n$ be a quintasymptotic sequence in a
Noetherian ring $R$. Let $1\leq i\leq n$ and $\frak p\in
\overline{Q^\ast}((x_1,\dots, x_i))$. If $\frak q$ is a minimal
prime divisor of $(\frak p,(x_{i+1},\dots, x_n))$, then $\frak
q\in\overline{Q^\ast}((\mathbf{x}))$.
\end{cor}
\proof By Proposition 2.14, we have
$\overline{Q^\ast}((x_1,\dots, x_j))=\overline{A^\ast}((x_1,\dots, x_j))$,
for all $1\leq j\leq n$. Therefore $x_{i+1},\dots,x_n$ is a
asymptotic sequence over $x_1,\dots, x_i$,  and so the conclusion follows from {\cite[Lemma 6.13] {McA1}}.\qed \\

\section{Quintasymptotic cograde}
In this section we show that the quintasymptotic cograde of an ideal
in a Noetherian local ring is unambiguously defined and behaves well
when passing to certain related local rings. We begin the following
useful lemma which is proved by L. J. Ratliff, Jr.

\begin{lem}
Let $(R, \frak m)$ be a complete Noetherian local ring that has only
one prime divisor of zero and let $I$ be an ideal of $R$. Then
$\ecogd(I)=\dim R/I$ and $\egd(I)+\ecogd(I)=\dim R$.
\end{lem}
\proof See \cite[Lemma 3.1]{Ra1}.

\begin{thm}
Let $I$ denote an ideal in a Noetherian local ring $(R, \frak m)$.
Then any two maximal quintasymptotic sequences over $I$ have the
same length. In fact,

\begin{align*}
\qacogd(I)=&\mathrm{min}\{\dim R^\ast/{IR^\ast +z}\mid z\in \mAss_{R^\ast} R^\ast\}\\
=&\mathrm{min}\{\dim R^\ast/z-\height(IR^\ast+z/z)\mid z\in
\mAss_{R^\ast} R^\ast\}.
\end{align*}
\end{thm}
\proof Let $\mathbf{x}=x_1, \dots, x_n$  be a maximal
quintasymptotic sequence over $I$. Then by Proposition 2.8,
$\mathbf{x}$ is a maximal quintasymptotic sequence over $IR^\ast$.
Also by Proposition 2.6, the image of $\mathbf{x}$ in $R^\ast/z$ a
quintasymptotic sequence over ${IR^\ast+z}/z$, for all $z\in
\mAss_{R^\ast} R^\ast$ and this image is a maximal quintasymptotic
sequence over ${IR^\ast+z}/z$ for some such $z$. Now, as $R^\ast/z$
is a complete local domain, it follows that the image of
$\mathbf{x}$ in $R^\ast/z$ is a quintessential sequence over
${IR^\ast+z}/z$, and so for all $z\in \mAss_{R^\ast} R^\ast$, we
have $n\leq \qecogd({IR^\ast+z}/z)$; and the equality holds for some
such $z$. Therefore
$$n={\rm min}\{\qecogd(IR^\ast+z/z)\mid z\in \mAss_{R^\ast} R^\ast\}.$$
On the other hand, in view of Lemma 3.1,
$$\qecogd({IR^\ast+z}/z)=\dim R^\ast/{IR^\ast+z}.$$ Consequently it
follows that $\qacogd(I)$ is unambiguously defined and
$$\qacogd(I)={\rm min}\{\dim IR^\ast+z/z \mid z\in \mAss_{R^\ast} R^\ast\}.$$
Finally, the last equality follows from the fact that:
$$\dim R^\ast/{IR^\ast +z}=\dim (R^\ast/z)/({IR^\ast +z}/z)$$ and
$$\height({IR^\ast +z}/z)+\dim (R^\ast/z)/({IR^\ast +z}/z)=\dim
R^\ast/z.$$ \qed \\
The following theorem shows that $\qacogd(I)$ behaves nicely when
passing to certain related rings and ideals.

\begin{thm}
Let $I$ and $J$ be ideals in a Noetherian local ring $(R,\frak m)$. Then

$\rm (i)$ If $I\subseteq J$, then $\qacogd(I)\leq \qacogd (J)$.

$\rm (ii)$ If $\Rad(I)=\Rad(J)$, then $\qacogd(I)=\qacogd(J)$.

$\rm (iii)$ $\qacogd(I)={\rm min}\{\qacogd(I+z/z) \mid z\in \mAss_R
R\}$.

$\rm (iv)$ If $T$ is a faithfully flat Noetherian extension of $R$,
then $\qacogd(I)\leq \qacogd(IT_\frak n)$ for all prime ideals
$\frak n$ in $T$ lying over $\frak m$ and equality holds if $\height
\frak m=\height\frak n$. In particular $\qacogd(I)=\qacogd(IR^\ast)$

$\rm (v)$ If $T$ is a finite extension ring of $R$ such that every
minimal prime of $T$ lies over a minimal prime of $R$, then
$\qacogd(I)\leq\qacogd(IT_\frak n)$ for all prime ideals $\frak n$
in $T$ lying over $\frak m$ and equality holds for some such $\frak
n$.

$\rm (vi)$ $\qacogd(I)=\qacogd((I,X)R[X]_{(\frak m ,X)}).$

\proof (i) and (ii) follow from Theorem 3.2 and (iii) follows from
Theorem 3.2 and Proposition 2.6. In order to prove (iv), let
$\mathbf{x}=x_1, \dots, x_n$ be a maximal quintasymptotic sequence
over $I$ and let $\frak n$ be a prime ideal in $T$ lying over $\frak
m$. As $(I,(\mathbf{x}))T\subseteq \frak n$, it follows that
$\mathbf{x}$ is a quintasymptotic sequence over $IT_\frak n$, by
Propositions 2.5 and 2.8. Whence $\qacogd(I)\leq\qacogd(IT_\frak
n)$, by Theorem 3.2. Now, if $\height\frak m=\height\frak n$, then
$\frak n$ is a minimal prime over $\frak mT$. Thus $\frak nT_{\frak
n}\in \overline{Q^\ast}((I,(\mathbf{x}))T_{\frak n})$ by Lemma
2.3(iv),(vii). Therefore $\qacogd$$(I)=\qacogd(IT_\frak n)$, by
Theorem 3.2.

For $\rm (v)$, let $\mathbf{x}=x_1, \dots, x_n$ be a maximal
quintasymptotic sequence over $I$ and  let $\frak n$ be a prime ideal in
$T$ lying over $\frak m$. Since $\frak n$ contains
$(I,(\mathbf{x}))T$, it follows that $\mathbf{x}$ is a
quintasymptotic sequence over $IT_\frak n$, by Propositions 2.5 and
2.9. Whence $\qacogd(I)\leq\qacogd(IT_\frak n)$ for all prime ideals
$\frak n$ in $T$ lying over $\frak m$, by Theorem 3.2. Also, by
Proposition 2.9, there exists a maximal ideal $\frak n'$ in $T$ such
that $\frak n'\in \overline{Q^\ast}((I,(\mathbf{x}))T)$. Thus $\frak
n'T_{\frak n'}\in \overline{Q^\ast}((I,(\mathbf{x}))T_{\frak n'})$
by Lemma 2.3(iv), and hence $\mathbf{x}$ is a maximal
quintasymptotic sequence over $IT_{\frak
n'}$. Therefore $\qacogd(I)=\qacogd(IT_{\frak n'})$ by Theorem 3.2. Finally, (v) follows immediately from Proposition 2.10 and Theorem 3.2. \qed \\
\end{thm}

In remainder of this section we give several bounds on
quintasymptotic cograde of an ideal.

\begin{thm}
Let $I$ be an ideal in a Noetherian local ring $R$. Then the
following hold:

$\rm (i)$ $\qacogd(I)\leq\dim R/I$.

$\rm (ii)$ $\qacogd(I)\geq {\rm min}\{\dim R^\ast/\frak
q\mid \frak q\in\overline{Q^\ast}(IR^\ast)\}$.
\end{thm}
\proof For (i) we have,
\begin{center}
$\dim R/I=\dim R^\ast/IR^\ast\geq \min\{\dim R^{\ast}/{IR^{\ast}+z}
\mid z\in \mAss_{R^\ast} R^\ast\},$
\end{center}
and so, by Theorem 3.2, we have $\qacogd(I)\leq\dim R/I$.

In order to prove (ii), in view of Theorem 3.2, there exists a minimal prime $z$ of $R^\ast$
such that $\qacogd(I)=\dim R^\ast/{IR^\ast+z}$. Thus there exists a
minimal prime divisor $\frak q$ of $IR^\ast+z$ such that $\qacogd(I)=\dim R^\ast/\frak q$.
Hence, by Lemma 2.3(vi), $\frak q\in \overline{Q^\ast}(IR^\ast)$, and so the result follows.\qed \\

\begin{lem}
Let $I$ and $J$ be ideals in a Noetherian ring $R$ and let $\frak
p\in\Spec R$ such that $I\subseteq J\subseteq \frak p$ and $\frak
p\in\overline{Q^\ast}(I)$.  Then $\frak p\in\overline{Q^\ast}(J)$.
\end{lem}
\proof  It follows from Definition 2.1.\qed\\

\begin{prop}
Let $I$ be an ideal in a Noetherian local ring $(R,\frak m)$ and let
$y_1,\dots, y_k$ be an asymptotic sequence of elements of $R$ such
that $y_j\in I$ for all $j(1\leq j\leq k)$. Then there is a maximal
quintasymptotic sequence over $I$, say $ x_1,\dots, x_n$, such that
$y_1,\dots, y_k, x_1,\dots, x_n$ is an asymptotic sequence in $R$.
In particular $\qacogd(I)+\qagd(I)\leq \qagd(\frak m)$.
\end{prop}
\proof Let $n$ be the length of a maximal quintasymptotic sequence
over $I$. If $n=0$ we are done. If $n>0$, then $\frak m \notin
\overline{Q^\ast}(I)$, and so by Lemma 3.5, $\frak m \notin
\overline{Q^\ast}((y_1,\dots, y_k))$. Pick $x_1\in \frak m$ with
$x_1\notin \cup \{ \frak p \mid \frak p\in \overline{Q^\ast}(I)\}$
and $x_1\notin \cup \{\frak p \mid \frak p\in
\overline{Q^\ast}((y_1,\dots, y_k))\}$. Now, $x_1$ is a
quintasymptotic sequence over $I$ and the length of a maximal
quintasymptotic sequence over $(I,(x_1))$ is $n-1$. Since the choice
of $x_1$ assure that $y_1,\dots, y_k,x_1$ is an asymptotic sequence
in $R$,
we now may use induction.\qed \\

The next result determines when the inequality in Proposition 3.6 is
equality.

\begin{thm}
Let $I$ be an ideal in a Noetherian local ring $(R,\frak m)$. Then
the following statements are equivalent:

$\rm (i)$ $\qacogd(I)+\agd(I)=\agd(\frak m)$.

$\rm (ii)$ For all $z\in \mAss_{R^\ast} R^\ast$ with $\agd(\frak
m)=\dim R^\ast/z$, we have $\agd(I)=\height({IR^\ast+z}/z)$ and
$\qacogd(I)=\dim R^\ast/z-\height({IR^\ast+z}/z)$.

$\rm (iii)$ There exists $z\in \mAss_{R^\ast} R^\ast$ such that
$\agd(I)=\height({IR^\ast+z}/z)$ and $\qacogd(I)=\dim
R^\ast/z-\height({IR^\ast+z}/z)$.
\end{thm}

\proof For prove $\rm (i)\Longrightarrow \rm (ii)$, let
$\qacogd(I)+\agd(I)=\agd(\frak m)$ and let $z\in \mAss_{R^\ast}
R^\ast$ such that $\agd(\frak m)=\dim R^\ast/z$. Then  by Theorem
3.2 and {\cite[Proposition 6.10] {McA1}}, we have
$$\qacogd(I)\leq\dim R^\ast/z-\height({IR^\ast+z}/z)\leq \agd(\frak
m)-\agd(I)=\qacogd(I).$$ Therefore, $\agd(I)=\height({IR^\ast
+z}/z)$ and $\qacogd(I)=\dim R^\ast /z-\height({IR^\ast +z}/z)$.

In order to prove the implication $\rm (ii)\Longrightarrow \rm
(iii)$, it suffices to show that, there exists $z\in \mAss_{R^\ast}
R^\ast$ such that $\agd(\frak m)=\dim R^\ast /z$. To this end we use
{\cite[Proposition 6.10] {McA1}}.
Finally, the implication  $\rm (iii)\Longrightarrow \rm (i)$ is obvious, by {\cite[Proposition 6.10] {McA1}} and Proposition 3.6. \qed \\

Now we give some lower bounds on $\qacogd(I)$.

\begin{rem}
Let $I$ be an ideal in a Noetherian ring $R$. Then the following
hold:

{\rm(i)} It follows from the Lemma {\rm 2.3(iii)}, Corollary {\rm2.15} and
{\cite[Proposition 3.10] {KR}}, that $$\egd(I)=\qeg(I)\leq \qagd(I)=\agd(I).$$

{\rm(ii)} If $R$ is local, then by the Lemma {\rm 2.3(iii)}, we have,
\begin{center}
$\ecogd(I)\leq\qecogd(I)\leq \qacogd(I)$ and $\ecogd(I)\leq \acogd(I)\leq
\qacogd(I).$
\end{center}
\end{rem}

\begin{cor}
Let $I$ be an ideal in a Noetherian local ring $(R,\frak m)$. Then, for all large $n$,
$\qacogd(I)\geq \egd(\frak m/I^n)$.
\end{cor}
\proof This is clear by {\cite[Theorem 5.5] {Ra1}} and the Remark 3.8. \qed \\

\begin{cor}
Let $I$ be an ideal in a Noetherian local ring $(R,\frak m)$. Then, for all large $n$,
$\qacogd(I)\geq \grade(\frak m/I^n)$.
\end{cor}
\proof This is clear by Corollary 3.9 and  Lemma {\rm 2.3(iii)}.\qed \\

\begin{cor}
Let $I$ be an ideal in a Noetherian local ring $(R,\frak m)$ such
that for all large $n$, $(R/I^n)^\ast$ has no imbedded prime
divisors zero. Then  $\qacogd(I)\geq \agd(\frak m/I^n)$,
for all large $n$.
\end{cor}
\proof This is clear by Corollaries 3.9 and 2.16.\qed \\

This section will be closed with the another lower bound on
$\qacogd(I)$ in connection with analytic spread. Let us, firstly,
recall the important notion of the {\it analytic spread}
of $I$ in a local ring $(R, \frak{m}$), which introduced by
Northcott and Rees in \cite{NR} and defined as $l(I):= \dim\,\mathscr{R}/ (\frak{m},u)\mathscr{R}.$

\begin{thm}
Let $I$ be an ideal in a Noetherian local ring $(R,\frak m)$. Then,
$$\qacogd(I)\geq \agd(\frak m)-l(I).$$
\end{thm}
\proof By Theorem 3.2 there exists $z\in \mAss_{R^\ast} R^\ast$ such
that $$\qacogd(I)=\height({\frak mR^\ast}/z)-\height{IR^\ast+z}/z.$$
Also, by {\cite[Proposition 6.10] {McA1}}, we have $\agd(\frak
m)\leq \height({\frak mR^\ast}/z)$. In other hand, by {\cite[Remark
5.5.4] {KR}} and {\cite[Lemma 6.5] {McA1}},
$$\height({IR^\ast+z}/z)\leq l({IR^\ast+z}/z)\leq l(I).$$
Now the desired result follows. \qed \\

\section{Quintasymptotic cograde and unmixedess}

In this section we use quintasymptotic cograde to obtain some characterizations of quasi-unmixed rings and
another related class of local rings. We begin with the quasi-unmixed case.

\begin{prop}
Let $(R,\frak m)$ be a Noetherian local ring. Then the following are
equivalent:

$\rm (i)$ $R$ is quasi-unmixed.

$\rm (ii)$ $\agd(I)=$ $\height I$ for all ideals $I$ of $R$.

$\rm (iii)$ $\agd(\frak m)=\height\frak m$.

$\rm (iv)$ $\qacogd(0)=\dim R$.

$\rm (v)$ $\qacogd(I)=\dim R/I$ for all ideals $I$ of $R$
generated by an asymptotic sequence of elements of $R$.
\end{prop}
\proof Since by {\cite[Lemma 2.5] {Rat1}} every quasi-unmixed local
ring is locally quasi-unmixed ring, it follows from {\cite[Corollary
5.8] {McA1}} that (i)-(iii) are equivalent. Assume that (ii) holds
and let $\mathbf{x}=x_1, \dots, x_n$ be an asymptotic sequence in
$R$ and let $I:=(\mathbf{x})$. Then there are elements
$x_{n+1},\dots,x_r$ in $R$ such that
$x_1,\dots,x_n,x_{n+1},\dots,x_r$ is a maximal asymptotic sequence
in $R$, and so $\qacogd(I)=r-n$. By assumption and Corollary 2.4(i),
we have $\height I=n$ and $\height\frak m=r$. Hence $n+\dim R/I\leq
\dim R=\height\frak m=r$, and so $\dim R/I\leq r-n=\qacogd(I)$. Now,
the implication $\rm (ii)\Longrightarrow\rm (v)$ follows from
Theorem 3.4. The implications
$\rm (v)\Longrightarrow\rm (iv)$ and $\rm (iv)\Longrightarrow\rm (iii)$ obviously are true. \qed \\

\begin{thm}
Let $(R, \frak m)$ be a Noetherian local ring. Consider the following
conditions:

$\rm (i)$ $\qacogd(I)=\dim R/I$ for every ideal $I$ of $R$.

$\rm (ii)$ $\agd(I)+\qacogd(I)=\dim R$ for every ideal $I$ of $R$.

$\rm (iii)$ $\agd(\frak p)+\qacogd(\frak p)=\dim R$
for every prime ideal $\frak p$ of $R$ with $\dim R/\frak p=1$.

$\rm (iv)$ $R$ is quasi-unmixed and there exists a unique minimal
prime divisor in $R$.\\
Then $\rm (i)\Longleftrightarrow\rm (ii)\Longleftrightarrow\rm
(iii)\Longrightarrow\rm (iv)$.
\end{thm}
\proof  $\rm (i)\Longrightarrow\rm (iv)$: If (i) holds and $I=(0)$, then by Proposition 4.1, $R$ is a
quasi-unmixed ring. Let $\frak p$ and $\frak q$ be distinct minimal
prime divisors in $R$. Then, by {\cite[34.5] {Nag1}}, we have
$$\dim R/\frak p=\dim R/\frak q=\dim R.$$
Let $z$ be a minimal prime divisor of $\frak qR^\ast$ such that
$\dim R^\ast/z=\dim R^\ast/\frak qR^\ast$, and so $\dim
R^\ast/z=\dim R$. Thus $z\in \mAss_{R^\ast} R^\ast$. Therefore, by
Theorem 3.2, we have
$$\qacogd(\frak p)\leq\dim R^\ast/{\frak pR^\ast+z}.$$

As $\frak p$ and $\frak q$ are distinct, it yields that
$$\dim R^\ast/{\frak pR^\ast+z}<\dim R^\ast/\frak qR^\ast=\dim R/\frak q=\dim R/\frak
p.$$
Hence $\qacogd (\frak p)<\dim R/\frak p$, which is a contradiction.

$\rm (i)\Longrightarrow\rm (ii)$: If (i) holds, then $R$ is a quasi-unmixed ring,
and so $\agd(I)=\height I$, by Proposition 4.1.  Thus $\rm(ii)$ is true by
{\cite[34.5] {Nag1}}. It is clear that $\rm (ii)\Longrightarrow\rm
(iii)$.

In order to prove the implication $\rm (iii)\Longrightarrow\rm (i)$,
suppose, the contrary, that $\rm(i)$ is not true. Then, there is an
ideal $I$ of $R$ such that $\qacogd (I)\neq\dim R/I$ and $\dim
R/I\leq\dim R/J$ for every ideal $J$ of $R$ with $\qacogd
(J)\neq\dim R/J$. Now, suppose $\dim R/I=d$. Then by Theorem 3.4(i),
$d>0$. If $d=1$, then $\qacogd(I)=0$ by Theorem 3.4(i), and so by
Theorem 3.3(i), $\qacogd(\frak p)=0$, for all prime ideals $\frak p$
in $R$ containing $I$, which is a contradiction. Therefore $d>1$ and
by {\cite[Proposition 2.2] {Rat2}} there exists infinitely many
prime ideals $\frak p$ in $R$ such that $I\subseteq \frak p$ and
$\dim R/\frak p=d-1$. Let $\mathscr{P}$ be the set of these prime
ideals. Then, by Theorem 3.2 there exists $z\in \mAss_{R^\ast}
R^\ast$ such that $\qacogd(I)=\dim R^\ast/{IR^\ast +z}$. Let
$$\mathscr{Q}=\{\frak q\in\Spec R^\ast\mid\text{there is}\,\frak p\in
\mathscr{P}\,\text{with}\,\frak pR^\ast+z\subseteq \frak q \, \text{and}
\dim R^\ast/{\frak pR^\ast +z}=\dim R^\ast/\frak q\}.$$  We now
show that $\mathscr{Q}$ is infinite and
$\dim R^\ast/\frak q=d-1$ for all $\frak q\in \mathscr{Q}$. To do this, let
$\frak q\in \mathscr{Q}$.  Then there exists $\frak p\in \mathscr{P}$
such that $\frak pR^\ast+z\subseteq \frak q$ and $\dim R^\ast/{\frak
pR^\ast +z}=\dim R^\ast/\frak q$. By choice of $d$ we have
$\qacogd(\frak p)=\dim R/\frak p=d-1$. Therefore
$$d-1=\qacogd(\frak p)\leq \dim R^\ast/{\frak pR^\ast +z}=\dim R^\ast/\frak q.$$ In other hand, we have
$$\dim R^\ast/\frak q\leq\dim R^\ast/\frak pR^\ast=\dim R/\frak p=d-1.$$ Therefore
$\dim R^\ast/\frak q=d-1$ and $\frak q$ is a minimal prime divisor
of $\frak pR^\ast$. Consequently,  there are infinitely many $\frak q$,
since there are infinitely many $\frak p$. As
$$d-1=\dim R^\ast/\frak q\leq\dim R^\ast/{\frak pR^\ast +z}\leq\dim R^\ast/{IR^\ast+z}=\qacogd(I)<\dim R/I=d,$$
it follows that $\dim R^\ast/{IR^\ast+z}=d-1$. Hence $\frak q$ is a minimal prime
divisor of $IR^\ast+z$, and therefore $IR^\ast+z$ has infinitely many
minimal prime divisors, which is a contradiction. That is $\rm(i)$ holds. \qed

\begin{thm}

Let $R$ be a complete Noetherian local ring. Then the following conditions are
equivalent:

$\rm (i)$ $\qacogd(I)=\dim R/I$, for every ideal $I$ of $R$.

$\rm (i)$ $\agd(I)+\qacogd(I)=\dim R$, for every ideal
$I$ of $R$.

$\rm (iii)$ $\agd(\frak p)+\qacogd(\frak p)=\dim R$,
for every prime ideal $\frak p$ of $R$ with $\dim R/\frak p=1$.

$\rm (iv)$ $R$ has a unique minimal prime divisor of zero.
\end{thm}
\proof In view of Theorem 4.2, it suffices to show that $\rm (iv)\Longrightarrow\rm
(i)$; and this follows from Theorem 3.2. \qed \\

\begin{center}
{\bf Acknowledgments}
\end{center}

The authors would like to thank Dr. Monireh Sedghi for
her careful reading of the first draft and many helpful
suggestions. Also, we would like to thank from the Institute for
Studies in Theoretical Physics and Mathematics (IPM) for its
financial support.


\end{document}